\documentclass[11pt]{article}
\usepackage{amsmath}
\usepackage{amssymb}
\usepackage{amsthm}
\usepackage[margin=3.5cm]{geometry}
\usepackage{latexsym,amsfonts}
\usepackage{algorithm,algorithmic}
\usepackage{graphicx}
\usepackage{subfigure}
\usepackage{epstopdf}
\usepackage{multirow}

\def\fatn{\mathbf{n}}
\def\fatx{\mathbf{x}}

\def\fatnu{\mathbf{\nu}}
\def\voxel{\mathcal{V}}

\def\opM{\mathcal{M}}
\def\opD{\mathcal{D}}

\title{Local error estimates for adaptive simulation of the Reaction--Diffusion Master Equation via operator splitting}
\author{{\small 
      Andreas Hellander$^{\mbox{\tiny 1,}}$\footnote{To whom correspondence should be addressed. email: \tt andreash@cs.ucsb.edu}, Michael Lawson$^{\mbox{\tiny 1}}$, Brian Drawert$^{\mbox{\tiny 1}}$, Linda Petzold$^{\mbox{\tiny 1}}$  
       }
}

\begin{document}
\maketitle
\vspace{-20pt}

\begin{center}
{\footnotesize\em 
$^{\mbox{\tiny\rm 1}}$Department of Computer Science,
University of California Santa Barbara, Santa Barbara, CA 93106-5070, USA \\
[3pt]}
\end{center}

\begin{abstract}

The efficiency of exact simulation methods for the reaction-diffusion master equation (RDME) is severely limited by the large number of diffusion events if the mesh is fine or if diffusion constants are large. Furthermore,  inherent properties of exact kinetic-Monte Carlo simulation methods limit the efficiency of parallel implementations. Several approximate and hybrid methods have appeared that enable more efficient simulation of the RDME. A common feature to most of them is that they rely on splitting the system into its reaction and diffusion parts and updating them sequentially over a discrete timestep. This use of operator splitting enables more efficient simulation but it comes at the price of a temporal discretization error that depends on the size of the timestep. So far, existing methods have not attempted to estimate or control this error in a systematic manner. This makes the solvers hard to use for practitioners since they must guess an appropriate timestep. It also makes the solvers potentially less efficient than if the timesteps are adapted to control the error. Here, we derive estimates of the local error and propose a strategy to adaptively select the timestep when the RDME is simulated via a first order operator splitting. While the strategy is general and applicable to a wide range of approximate and hybrid methods, we exemplify it here by extending a previously published approximate method, the Diffusive Finite-State Projection (DFSP) method, to incorporate temporal adaptivity. 
\end{abstract}

\section{Introduction}

To understand biological systems on the cellular level, it is often essential to account for the impact of noise due to small molecule count. 
This was beautifully demonstrated in the now classic results on the effect of stochasticity in gene regulatory systems \cite{StochasticGene97, Elowitz16082002}.  Spatial distribution of molecules in a cell can result in locally small populations of key chemical species, such that noise drives essential behavior, as in the case of symmetry breaking across many eukaryotic cell types \cite{under_det}. 
Spatial stochastic modeling has already begun to yield new insights in systems such as spatiotemporal oscillators \cite{FaEl,Sturrock1:2013,Sturrock2:2013}, MAPK signaling \cite{TaTNWo10}, self-organization of proteins into clusters \cite{BISTAB} and polarization of proteins on the cell membrane in yeast \cite{Altschuler:2008}. 

On a mesoscopic level, the reaction-diffusion master equation (RDME) is frequently used to model these systems. 
Space is subdivided into subvolumes that can individually be treated as well-mixed. Reactions within a subvolume are expressed in the form of the chemical master equation (CME) \cite{GillMaster} and realizations of the process can be generated using Gillespie's stochastic simulation algorithm (SSA) \cite{SSA}. 
Molecules can move freely between neighboring voxels via diffusive jumps, which are modeled as linear jump events in a Markov process. 
Optimized exact simulation methods such as the Next Subvolume Method (NSM) \cite{BISTAB} can be used to generate statistically correct realizations of the RDME. 

As with all exact methods applied to RDME models, the NSM suffers from a potentially high computational cost due to having to explicitly simulate each diffusion event. The number of diffusive transfers between voxels grows rapidly as the mesh resolution is made finer, and as a result the majority of computation time tends to be spent on sampling diffusion events. Additionally, these methods are also inherently serial, which has thwarted attempts to increase efficiency via parallelization. 

To speed up simulation of the RDME many methods rely on operator splitting. 
By splitting the operators, most often with a Lie-Trotter scheme \cite{lts}, the reaction and diffusion steps can be solved independently. While diffusion carries the bulk of the computational burden in exact solvers, it is often possible to take advantage of the structure and linear nature of the discretized diffusion equation to speed up this step in an operator-split solver. 
Examples of approximate methods that have been proposed to speed up the simulation of the RDME by reducing the cost of the diffusive step include methods based on tau-leaping \cite{RosBayKou,marquez-lago:104101}, the multinomial simulation algorithm \cite{LaGiPe}, spatially adaptive hybrid methods \cite{FeHeLo2009} and the diffusive finite state projection method (DFSP) \cite{dfsp}. While splitting in itself does not resolve the issue of the inherent stiffness of the diffusion operator, the continued introduction of methods for simulating the RDME via operator splitting highlights 
the potential of this approach. Another recent use of operator splitting in the RDME context is the use of Lie-Trotter splitting to simulate fractional diffusion \cite{Bayati:2013}. Yet another advantage of splitting is that it converts a largely serial problem, which is known to be difficult to parallelize in an efficient manner, into a naturally parallelizable one.  For the existing approximate and hybrid methods for the RDME, splitting the physics (reaction and diffusion) is necessary. For parallel implementations, another possibility is to split the computational grid into blocks, as proposed by \cite{Arampatziz:2012}, where the error introduced by operator splitting at block boundaries was analyzed. Our analysis here is different since it applies to the case of splitting the reaction and diffusion operators.

Splitting the operators introduces an error that depends on the size of the splitting time step, however previous algorithms that rely on operator splitting have not attempted an {\it a priori} error estimator.  Without such an estimator these methods have no way to automatically control the splitting error. 
This limits their usefulness for practitioners. 
From an efficiency point of view, not knowing and controlling the error might lead to the use of unnecessary small timesteps at the price of slower simulations.  
In this work we seek to address theses issues
by presenting estimators of the local error in probability, mean and variance for a first order Lie-Trotter splitting of the RDME. 
The estimators allow control of the splitting error for spatial stochastic simulation and enables any method based on operator splitting to be implemented adaptively. 

This paper is organized as follows: in Section 2 we introduce the RDME. 
In Section 3 we outline how to simulate the RDME using operator splitting. 
We derive our estimator for the local error in the PDF, mean and variance, and demonstrate the accuracy of the estimator for an example problem. 
Finally, in Section 4 we present how the local error estimates can be used to extend an approximate method for the RDME to incorporate temporal adaptivity.

\section{Spatial Stochastic Simulation using the RDME}

Given a system with $N_s$ chemical species $X_s$ reacting in a physical domain $\Omega$, discretize $\Omega$ into $N_v$ non-overlapping voxels $\voxel_i$, with volume $|\voxel_i|$, and let the state of the system be described by the $(N_v \times N_s)$ state matrix $\fatx$, where the element $x_{is}$ is the copy number of species $X_s$ in voxel $\voxel_i$. Let $\fatx_{i,\cdot}$ denote the $i$-th row of $\fatx$ and $\fatx_{\cdot,s}$ the s-th column. The reaction network consists of $N_r$ chemical reactions $r = 1 \ldots N_r$. The $(1\times N_s)$ stoichiometry vector $\fatn_{ir}$ describes the change in the state, $\tilde \fatx_{i,\cdot} = \fatx_{i,\cdot} + \fatn_{ir}$, when reaction $r$ occurs in voxel $i$ and the propensity function for that reaction is $a_{ir}(\fatx_{i,\cdot})$. 

Diffusion of species $X_s$ along the edge (2D) or face (3D) connecting voxels $\voxel_i$ and $\voxel_j$ is modeled as a linear jump event, or diffusive transfer,  

\begin{equation}
X_{is} \xrightarrow{d_{ijs}} X_{js},
\label{eq:linjump}
\end{equation}
\noindent
with propensity function $\mu(\fatx) = d_{ijs}x_{is}$.  The change in state due to the diffusive transfer is described by the $(N_v\times 1)$ stoichiometry vector $\fatnu_{ijs}$, such that the new state is $ \tilde \fatx_{\cdot,s} = \fatx_{\cdot,s} + \fatnu_{ijs}$. All entries of $\fatnu_{ijk}$ are zero except $\fatnu_{ijs}(i)=-1$ and $\fatnu_{ijs}(j)=1$. In the case of reactions only, the probability density function $p(\fatx,t) \equiv p(\fatx,t | \fatx_0,0)$ obeys the master equation

\begin{align}
  \nonumber
  \frac{\mathrm{d}}{\mathrm{dt}}p(\fatx, t) &=
  \mathcal{M} p (\fatx, t) \equiv                                               \\
  \nonumber
  &\sum_{i=1}^{N_v}
  \sum_{r = 1}^{N_r}
  a_{ir}(\fatx_{i \cdot}-\fatn_{ir})
  p(\fatx_{1 \cdot},\ldots,\fatx_{i \cdot}-\fatn_{ir},
  \ldots,\fatx_{N \cdot}, t)                                              \\
  -&\sum_{i=1}^{N_v}
  \sum_{r = 1}^{N_r}
  a_{ir}(\fatx_{i \cdot})p(\fatx, t).
  \label{eq:RDMEre}
\end{align}
\noindent
For the case of one subvolume, $N_v=1$, \eqref{eq:RDMEre} reduces to the CME for a well stirred system. 
For a system with only diffusion, the master equation takes the form 

\begin{align}
  \nonumber
  \frac{\mathrm{d}}{\mathrm{dt}}p(\fatx, t) &=
  \mathcal{D} p (\fatx, t) \equiv                                               \\
  \nonumber
  &\sum_{s=1}^{N_s}
  \sum_{i = 1}^{N_v}
  \sum_{j=1}^{N_v}
  \mu(\fatx_{\cdot j}-\fatnu_{ijs})
  p(\fatx_{\cdot 1},\ldots,\fatx_{\cdot j}-\fatnu_{ijs},
  \ldots,\fatx_{\cdot N}, t)                                              \\
  -&\sum_{s=1}^{N_s}\sum_{i=1}^{N_v}
  \sum_{k = 1}^{N_v}
  \mu(\fatx_{\cdot j})p(\fatx, t),
  \label{eq:RDMEdf}
\end{align}

\noindent
For a system with both diffusion and reactions,  the full RDME is then given by
\begin{equation}
\frac{d}{dt}p(\fatx,t) = (\mathcal{M}+\mathcal{D})p(\fatx,t).
\label{eq:rdme}
\end{equation}

The values of the diffusion rate constants $d_{ijs}$ depend on the diffusion constant $\gamma_s$ of species $X_s$ and the shapes and sizes of voxels $\voxel_i$ and $\voxel_j$. 
Let $u_s(\zeta,t) = E[X_{s}/|\Omega|]$, i.e. the the expected value of the concentration process corresponding to \eqref{eq:RDMEdf}.
In the thermodynamic limit $X_{s}/|\Omega| \to \infty$, $u_s$ is governed by the  diffusion equation
\begin{equation}
\frac{\partial}{\partial t} u_s(\zeta,t) = \Delta u_s(\zeta,t). 
\label{eq:diffusionpde}
\end{equation}
\noindent
Here, $\zeta$ is a spatial coordinate in a Cartesian coordinate system. Similarly, the Fokker-Planck equation that describes the time evolution of the probability density of a single particle undergoing Brownian motion is  given by

\begin{equation}
\frac{\partial}{\partial t} p(\zeta,t) = \Delta p(\zeta,t). 
\label{eq:diffusionpde}
\end{equation}
\noindent

\begin{figure}[H]
\centering
\subfigure[]{\includegraphics[width=0.30\linewidth]{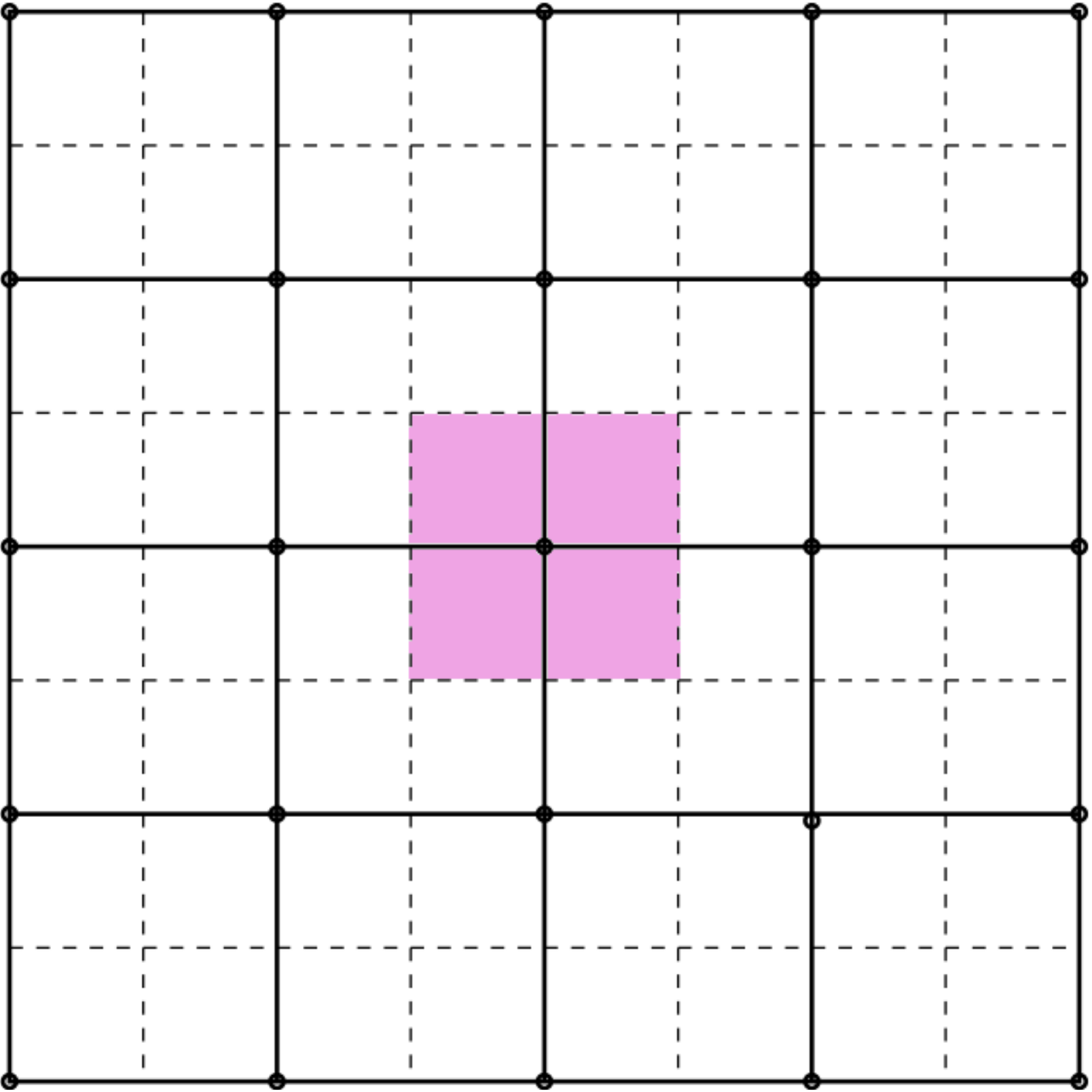}}
\hspace{0.5cm}
\subfigure[]{\includegraphics[width=0.30\linewidth]{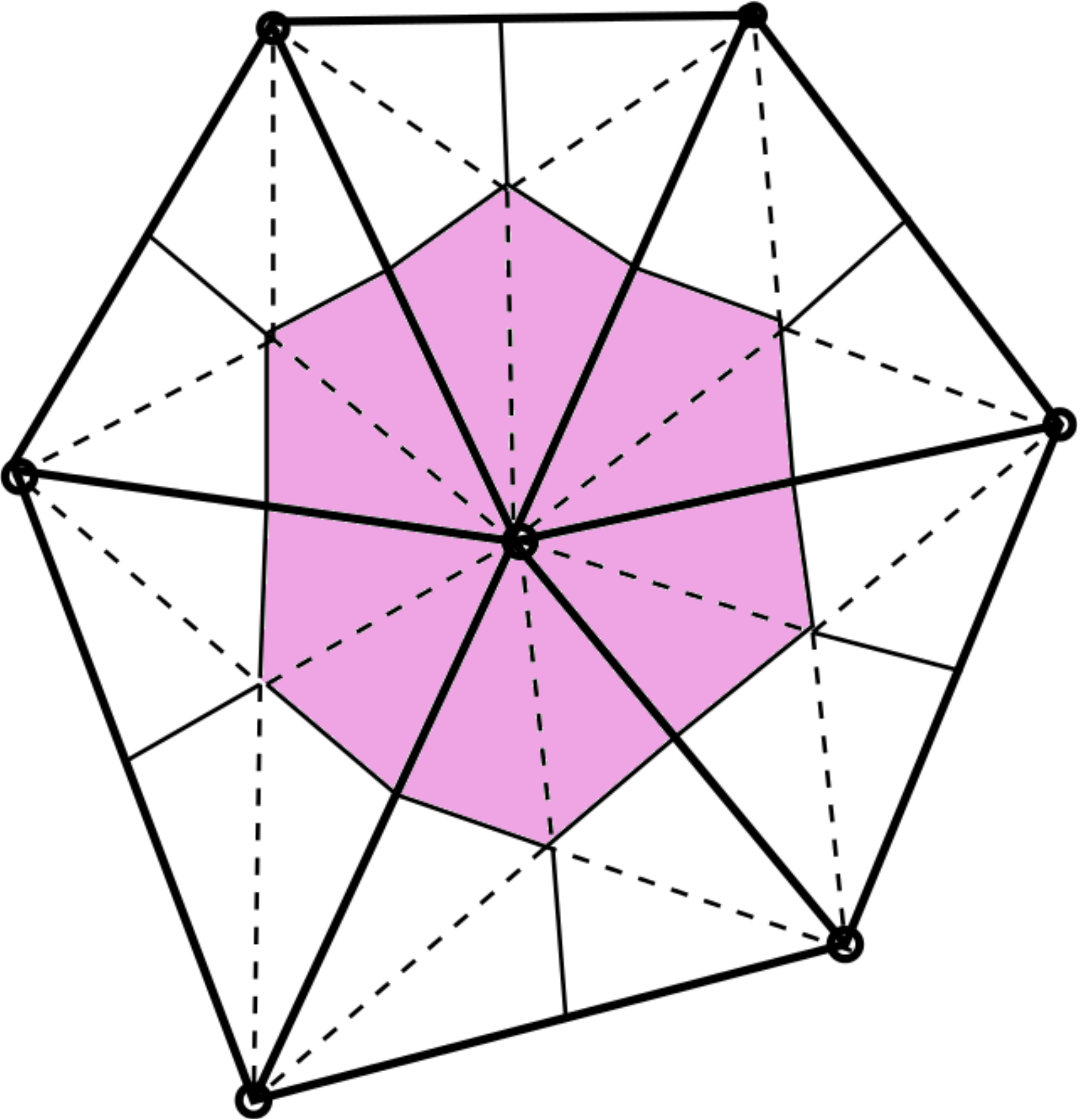}}
\caption{The molecules are assumed to be well mixed in the dual elements (depicted in pink) of the primal mesh  (solid thick lines). On a Cartesian grid (a), the duals are simply the volumes of the staggered grid. The dual of the triangular mesh in (b) is obtained by connecting the midpoints of the edges and the centroids of the triangles. In the conventional RDME, molecules are allowed to jump between immediate neighboring voxels.} 
\label{fig:meshes}
\end{figure}

\noindent 
Hence, by choosing $d_{ijs}$ according to a consistent spatial discretization of  the operator $\Delta$, we obtain mesoscopic jump coefficients that are motivated both from a microscale and macroscale perspective. To comply with the description of the RDME, the discretization used needs to couple only nearest neighbors in the mesh. For a uniform, Cartesian mesh such as shown in Figure \ref{fig:meshes} (a), the most natural choice is a centered finite difference discretization, giving $d_{ijs} = \gamma_s/h^2$ where $h$ is the side length of the voxels and $\gamma_s$ is the diffusion constant of species $s$. For unstructured, tetrahedral meshes which will be used later in this paper, $d_{ijs}$ can be obtained from a finite element (FE) discretization using linear Lagrange elements. For a detailed description of that methodology, see \cite{SPDEPEFHL}. Figure \ref{fig:meshes} (b) shows a triangular, unstructured mesh. The local volume where the molecules are assumed to be well mixed are given by the dual mesh.

\section{Operator splitting method for the RDME}

In this section we describe how an operator splitting method can be applied to generate realizations from the RDME with a controlled temporal discretization error. 

In the remainder of this paper we will tacitly assume that we can impose a bound on the state space. For every species in the system we assume that $P(X_{is}> x_{max},t)=0$, for all $t$. The state space is then finite (but very large). Technically, this will not necessarily hold true for an open system, but from a biophysical perspective it is a reasonable assumption since a system would require infinite energy to blow up. In the finite case, the operators in the RDME can be represented by finite matrices. 
Hereafter we will use the notation $\mathcal{M}$ and $\mathcal{D}$ to mean representations of the bounded operators resulting from the above truncation of the state space. The solution of \eqref{eq:rdme} can then simply be written

\begin{equation}
p(\fatx,t) = e^{t(\opM+\opD)}p(\fatx,0). 
\end{equation}
\noindent
We point out that even though this representation of the solution is simple, it is not feasible to solve \eqref{eq:rdme} this way since the state state space is too large except for trivial models and discretizations. In the well mixed case, however, deterministic methods to solve the CME have been developed for small to medium-sized models, see for example \cite{munsky:044104,BURRFSP,macnamara:1146,Jahnke2,mateescu:441}. 

Using a first order splitting method (Lie-Trotter splitting), $p(\fatx,t+\Delta t)$ can be approximated by
\begin{equation}
p_s(\fatx,t+\Delta t) = e^{\Delta t\opM}e^{\Delta t\opD}p(\fatx,t). 
\label{eq:lietrotter}
\end{equation} 
\noindent
Simulation based on \eqref{eq:lietrotter} proceeds in two half-steps, with the diffusion operator and the reaction operator acting sequentially
\begin{align}
\nonumber
&1.~~p_s^{n+1/2}=e^{\Delta t\opD}p_s^n(\fatx,t)\\
&2.~~p_s^{n+1}(\fatx,t+ \Delta t)=e^{\Delta t\opM}p_s^{n+1/2}.
\label{eq:rdmefracstep}
\end{align}
\noindent
Provided that the time step $\Delta t$ is sufficiently small, the local error in the scheme \eqref{eq:rdmefracstep} is $\mathcal{O}(\Delta t^2)$ and from standard theory for numerical solution of differential equations, the global error in the PDF is then proportional to $\Delta t$, i.e. 
\begin{equation}
\|p(\fatx,t)-p_s(\fatx,t)\| \leq C\Delta t.  
\end{equation}
\noindent
We note that it is not in general necessary for the operators to be bounded for the splitting method to result in a global error proportional to $\Delta t$. Jahnke shows in \cite{Jahnke:2010} that the global error for Strang splitting applied to the CME is second order under certain assumptions on the chemical reactions. Unfortunately, as Jahnke points out, while those conditions can be expected to hold for many systems, they are not easily verified in practice. Error bounds for exponential operator splitting have also been studied in other contexts in e.g. \cite{Jahnke:2000}. 

In the next subsection we derive estimates of the local error. In a subsequent section we will then illustrate how these estimates can be used to enable temporal adaptivity in an approximate method for the RDME.  

\subsection*{Local error}

Following standard theory for the analysis of scheme \eqref{eq:rdmefracstep}, the local error in the PDF (i.e. in probability) in the n+1-th timestep is given by
\begin{equation}
\mathbf{\epsilon}^{n+1} = (e^{\Delta t(\opM+\opD)}-e^{\Delta t\opM}e^{\Delta t\opD})p^n = \frac{\Delta t^2}{2}{[\opD,\opM]p^n} + \mathcal{O}(\Delta t^3).
\label{eq:eserrpdf}
\end{equation} 
where $[\opD,\opM]=(\opD\opM-\opM\opD)$ is the commutator of the operators. This follows directly from a series expansion of the left hand side in $\Delta t$.  Computing $\epsilon^{n+1}$ from \eqref{eq:eserrpdf} is obviously not tractable in general since knowledge of $p^n$ requires the solution of the full RDME. Even if an approximation of $p^n$ were available, the state space is still very large and the cost of evaluating the commutator would be prohibitive. Also, our goal is a strategy to estimate the error during the course of the generation of individual trajectories, whereas an estimation of $p^n$ would require very large ensembles of trajectories. By conditioning on the currently observed state $\fatx^n$ in timestep $n$ of a specific realization, we obtain the following sequence of conditional errors,  
\begin{equation}
\mathbf{\mathcal{E}}^{n+1} = (\mathbf{\epsilon}^{n+1} | X(t^n)=\fatx^n) = \frac{\Delta t^2}{2}{[\opD,\opM]\delta(\fatx^n)} + \mathcal{O}(\Delta t^3),
\label{eq:eserrpdfcond}
\end{equation}
\noindent
where $\delta$ is the Dirac delta function, hence $\delta(\fatx^n)$ corresponds to $P(X(t^n)=\fatx^n)=1$. Note that  $\mathbf{\mathcal{E}}^n$ is a random vector for each $n$ with unconditional expectation $E[\mathbf{\mathcal{E}^{n+1}}] = \epsilon^{n+1}$.

While there are many possible states that can have a non-zero value after one application of the commutator, hereafter called reachable states, it is possible to obtain simple and computable expressions for $\mathbf{\mathcal{E}}^{n+1}$ in \eqref{eq:eserrpdfcond}, due to the sparsity of $\delta(\fatx^n)$.  

The $L = 1+2N_eN_s+N_vN_r+2N_vN_rN_eN_s$  reachable states $\tilde{\fatx_l}$ in $[\opD,\opM]\delta(\fatx^n)$ are $\fatx^n$, $\fatx^n + \fatnu_{ijs}$, $\fatx^n+\fatn_{kr}$ and $\fatx^n+\fatn_{kr} + \fatnu_{ijs}$, $s=1\ldots N_s$, $r=1\ldots N_r$, $i,j,k = 1,\ldots, N_v$ and $N_e$ is the number of connections in the mesh. 
Below, the term to the left of the colon is a reachable state and the term to the right is 
the value of that state after applying the indicated operator to $\delta(\fatx^n)$. 

For example, the application of the diffusion operator $\opD$ on $\delta(\fatx^n)$, $\opD \delta(\fatx^n)$ results in 
\begin{align}
\fatx^n + \fatnu_{ijs} 	&: d_{ijs}x_{is}\quad \quad
 \label{opD_in}	\\
 \fatx^n & :  -d_0(\fatx^n)	,	
 \label{opD_out}		
\end{align}

\noindent
where $i,j$ are connected subvolumes and $d_0(\fatx^n) \equiv \sum_{s=1}^{N_s}\sum_{i=1}^{N_v}\sum_{j\neq i}^{N_v} d_{ijs}x_{is}.$ 
Equation (\ref{opD_in}) enumerates the new states with non-zero value after one application of $\opD$ to $\delta(\fatx^n)$. 

Applying $\opM$ to $\opD\delta(\fatx^n)$ then gives the following reachable states and associated values. 

\begin{align}
\fatx^n : & ~ a_0(\fatx^n)d_0(\fatx^n)\\ 
\fatx^n + \fatn_{kr} : &-a_{kr}(\fatx^n) d_0(\fatx^n),\\
 \fatx^n + \fatnu_{ijs} 	:& -a_0(\fatx^n+\fatnu_{ijs})d_{ijs}x_{is},\\
  \fatx^n + \fatnu_{ijs} +\fatn_{kr} : &~a_{kr}(\fatx^n+\fatnu_{ijs})d_{ijs}x_{is},
  \label{MD_xvn}
\end{align}
\noindent
where we have defined $a_0(\fatx^N) = \sum_{k=1}^{N_v}\sum_{r=1}^{N_r} a_{kr}(\fatx)$.  
By computing the analogous values for $\opD(\opM(\fatx^n)\delta(\fatx^n))$ and taking the difference, we obtain the following expressions for the commutator error $\mathbf{\mathcal{E}^{n+1}}$ for the reachable states $\tilde{\fatx}_l$

\begin{align}
\fatx^n+\fatn_{kr} + \fatnu_{ijs} &:   \frac{\Delta t^2}{2}\left\{ \begin{array}{ll} 
\label{eq:commrdme1}
d_{ijs}x_{is}(a_{ir}(\fatx_{i,\cdot})-
a_{ir}(\fatx_{i,\cdot}+\fatnu_{ijs}))+d_{ijs}n_{irs}a_{ir}(\fatx_{i,\cdot})& k=i\\
d_{ijs}x_{is}(a_{jr}(\fatx_{j,\cdot})-a_{jr}(\fatx_{j,\cdot}+\fatnu_{ijs}))& k=j\\ 
0 & k\neq \{i,j\}
\end{array}\right.\\ 
\fatx^n+\fatn_{ir} &: - \frac{\Delta t^2}{2}\sum_{s=1}^S\sum_{j\neq i}^N d_{ijs}n_{irs}a_{ir}(\fatx_{i,\cdot})\label{eq:commrdme2}\\ 
\fatx^n+\fatnu_{ijs}  &: \frac{\Delta t^2}{2} \sum_{r=1}^R d_{ijs}x_{is}\left(a_{ir}(\fatx_{i,\cdot}+\fatnu_{ijs})-a_{ir}(\fatx^n_{i,\cdot}) + a_{jr}(\fatx_{j,\cdot}+\fatnu_{ijs})-a_{jr}(x_{j,\cdot})\right)\label{eq:commrdme3}\\
\fatx^n & :  0.\label{eq:commrdmeend}
\end{align}

\noindent

\begin{figure}[htpb]
\centering
\includegraphics[width=\linewidth]{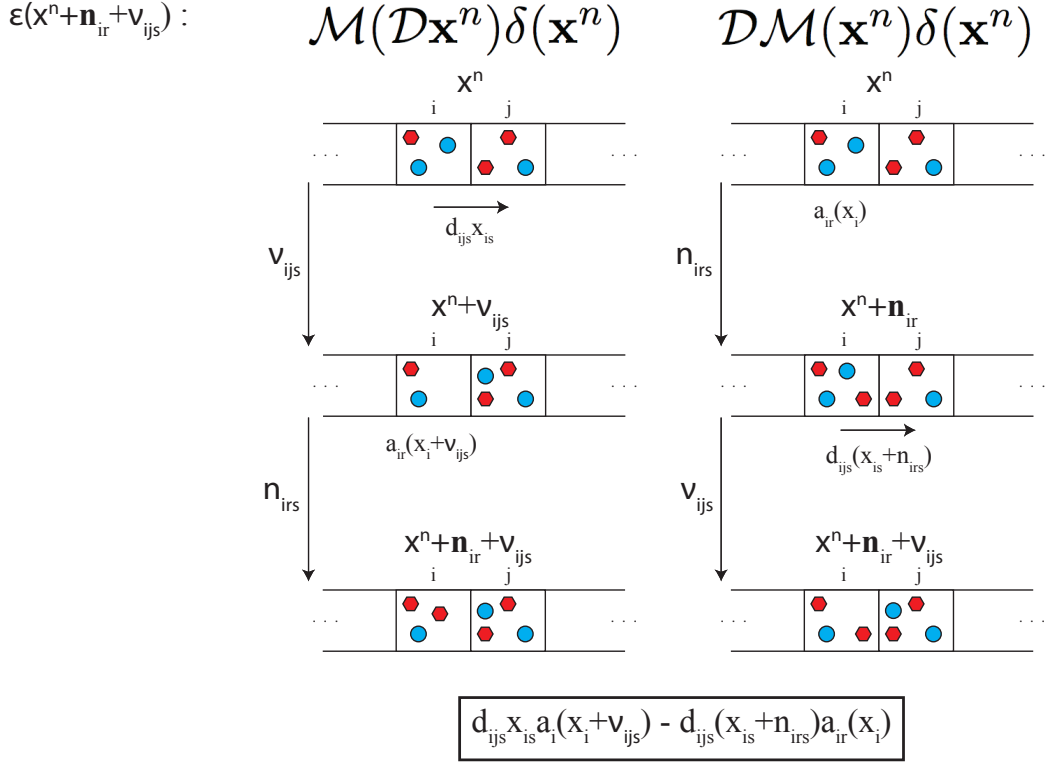}
\caption{
Illustration of the reachable state in Eq. (\ref{eq:commrdme1}) for $k=i$. 
On the left side, the diffusion operator is applied first, as in Eq. (\ref {opD_out}), followed by the reaction operator, resulting in the value derived in Eq. (\ref{MD_xvn}). 
On the right side the order of the operators is reversed, resulting in a change in the final value for the reachable state $\fatx^n+\fatn_{kr} + \fatnu_{ijs}$ (note that the reachable state is, by definition, the same in both orders, but the value associated with that state is different). 
The terms in the box represent the difference between the two orders of applying the operators and show terms that do not commute. 
Note that if we had defined the reachable state as $k\neq i,j$ in Eq. (\ref{eq:commrdme1}) then the reaction would occur in neither the originating subvolume nor the destination subvolume for the diffusion event and the terms would commute (and the value in the box would be zero).
} 
\label{fig:commutator}
\end{figure}

As seen in the expression for the cross-terms in \eqref{eq:commrdme1}, large parts of the operators commute.
Only terms pertaining to degrees of freedom that are sharing an edge in the mesh will be non-zero. Furthermore, \eqref{eq:commrdme1} will  be non-zero only if the reaction described by $\fatn_{kr}$ depends on species $s$, i.e. $n_{krs}$ is non-zero. A graphical representation  of the reachable states is given in Fig. \ref{fig:commutator}. 

\paragraph{Local error in mean and variance.}

Equations \eqref{eq:commrdme1} -- \eqref{eq:commrdmeend} give expressions for the local error in the PDF in one timestep (at $t^{n+1}$), given that we observe state $\fatx^n$ at time $t^n$. The obvious advantage of computing these expressions directly is that they can be used to derive many estimates of the error such as $l_1$, $l_\infty$ in PDF, Kolmogorov distance for marginal distributions, or moments. However, they are rather expensive to compute directly, even with appropriate optimizations. In many cases it is sufficient to control a weak error such as the errors in mean and variance. This is particularly true if the error in the individual methods for propagating the reaction operator and the diffusion operator is controlled only in a weak sense. For example, this would be the case if the diffusion operator were updated with $\tau$-leaping such as in \cite{FeHeLo2009}.

Based on \eqref{eq:commrdme1} -- \eqref{eq:commrdmeend} we can calculate the local error in the expected value of an arbitrary bounded function $g(X_{is})$ 

\begin{equation}
\Delta E[g(X_{is})] \equiv E[g(X_{is})]-E_{\tilde{f}}[g(X_{is})] =  \sum_{k=1}^Kg(x^k_{is})\mathcal{E}_k,
\end{equation}
where the subscript $\tilde{f}$ denotes expectation under the approximate PDF obtained from solving with the operator split method \eqref{eq:rdmefracstep}. For brevity, we have here dropped the superscript $n+1$ on $x$ and $\mathcal{E}$.

We introduce the notation $\Delta a_{irs}^{-}=a_{ir}(\fatx_{i,\cdot}+\fatnu_{ijs})-a_{ir}(\fatx_{i,\cdot})$ and $\Delta a_{irs}^{+}=a_{ir}(\fatx_{i,\cdot}+\fatnu_{jis})-a_{ir}(\fatx_{i,\cdot})$ and $\Delta g(y) = g(x_{is}^n+y)-g(x_{is}^n)$. By summing up outflow and inflow of differences in probability to state $X_{is}$ and using the fact that $\Delta E[g(X_{is})] = \Delta E[g(X_{is})-g(x_{is}^n)]$, we obtain
\begin{align}
\nonumber
(2\Delta t^{-2})&\Delta E[g(X_{is})] = \Delta g(-1)\sum_{j\neq i}^Nd_{ijs}x_{is}\sum_{r=1}^{N_r} \Delta a_{irs}^-+
\Delta g(+1)\sum_{j\neq i}^N d_{jis}\sum_{r=1}^{N_r} n_{jrs}a_{jr}(\fatx_{j,\cdot})\\
\nonumber
&+ \Delta g(+1)\sum_{j\neq i}^N d_{jis}x_{js}\sum_{r=1}^{N_r}\Delta a_{ris}^+
-\sum_{r=1}^{N_r} \Delta g(n_{irs})\sum_{j\neq i}^N d_{ijs}n_{irs}a_{ir}(\fatx_{i,\cdot})\\
\nonumber
&-\sum_{r=1}^{N_r}\Delta g(n_{irs})\sum_{j \neq i}\sum_{s' \neq s}d_{ijs'}x_{is'}\Delta a_{irs'}^- - 
\sum_{r=1}^{N_r} \Delta g(n_{irs})\sum_{j \neq i}\sum_{s' \neq s}d_{jis'}x_{js'}\Delta a_{irs'}^+\\
\nonumber
&+\sum_{r=1}^{N_r} \Delta g(n_{irs}-1)\sum_{j\neq i}^N (d_{ijs}n_{irs}a_{ir}(\fatx_{i,\cdot})- d_{ijs}x_{is}\Delta a_{irs}^-)\\
&-\sum_{r=1}^{N_r} \Delta g(n_{irs}+1)\sum_{j\neq i}^Nd_{jis}x_{js}\Delta a_{irs}^+. 
\label{eq:weakerror}
\end{align}

\noindent
To obtain the error in the mean, set $g(x) = x$. After simplification, we obtain 
\begin{align}
\nonumber
\Delta E^{n+1}[X_{is}| \fatx^n]&=0.5\Delta t^2(DR(\fatx^n))_{is}-0.5\Delta t^2\sum_{r=1}^{N_r} \sum_{s'=1}^S\sum_{j\neq i}^Nn_{irs}d_{ijs'}x_{is'}\Delta a_{irs'}^- \\
\nonumber
&-0.5\Delta ^2\sum_{r=1}^{N_r} \sum_{s'=1}^S \sum_{j \neq i}^Nn_{irs}d_{jis'}x_{js'}\Delta a_{irs'}^+ \\
\nonumber
= &\quad
0.5\Delta t^2(DR(\fatx))_{is} - \\
\nonumber
&- 0.5\Delta t^2\sum_{r=1}^{N_r} \sum_{s'=1}^Sn_{irs}(\Delta a_{irs'}^-x_{is'}\sum_{j\neq i}^Nd_{ijs'} + \Delta a^+_{irs'}\sum_{j \neq i}^Nd_{jis'}x_{js'}) \\ 
=&\quad 0.5\Delta t^2 \left((D_sR(\fatx))_{is} -  \sum_{r=1}^{N_r} n_{irs}\sum_{s'=1}^S\sigma_{irs'}\right).
\label{eq:splittingerrormean}
\end{align}
\noindent
where the operator $R(\fatx)$ is defined by
\begin{align}
R(\fatx)_{is} = \sum_{r=1}^{N_r}n_{irs}a_{ir}(\fatx_{i,\cdot})
\end{align}
\noindent
and
\begin{align}
\sigma_{irs'} = \Delta a_{irs'}^-x_{is'}\sum_{j\neq i}^Nd_{ijs'} + \Delta a^+_{irs'}\sum_{j \neq i}^Nd_{jis'}x_{js'}.
\end{align}
The $(N_v\times N_v)$ matrix $D_s$ is defined to have diagonal elements $d_{ii} = -\sum_{j\neq i} d_{ijs}$ and off-diagonal elements $d_{ij} = d_{ijs}$. 

Equation \eqref{eq:splittingerrormean} is an exact (up to $\mathcal{O}(\Delta t^3)$) formula for the estimate of the local error in mean for any functional form of the propensity functions $a_r(\fatx)$. For zeroth-order mass action reactions we have  $\sigma_{irs'} = 0$. For a first order reaction with a linear propensity (such as creation or monomolecular conversion) of the form $a_{ir}(\fatx_{i,\cdot}) = k_1x_{is}$ we have

\begin{align}
\nonumber
\sigma_{irs} = k_1(-x_{is}\sum_{j \neq i}d_{ijs} + \sum_{j \neq i} d_{jis}x_{js}) = k_1(D\fatx)_{is} = a_{ir}((D\fatx)_{i,\cdot}), n_{irs}\sigma_{irs'} = 0,
\end{align}
\noindent
and \eqref{eq:splittingerrormean} simplifies to 
\begin{align}
\Delta E^{n+1}[X_{is}] = 0.5\Delta t^2 \left( (DR(\fatx^n))_{is} -  R(D\fatx^n)_{is} \right),
\end{align}
\noindent
in agreement with the analogous expression for the commutator error for the reaction-diffusion PDE. 

To find the error in the second moment, set $g(x) = x^2$ in \eqref{eq:weakerror}, to obtain

\begin{align}
\nonumber
\Delta E[X_{is}^2] =
&\quad 0.5\Delta t^2(|D|R(x))_{is} + 2x_{is}\Delta E[X_{is}]  \\\nonumber
&-\Delta t^2\sum_r n_{irs}\sum_{j \neq i} (d_{ijs}x_{is}\Delta a_{irs}^- - d_{jis}x_{js}\Delta a_{irs}^+)\\ 
&- 0.5\Delta t^2\sum_{r} n_{irs}^2\sum_{s'=1}^S\sigma_{irs'} - \Delta t^2\sum_{j \neq i} d_{ijs}\sum_r n_{irs}^2a_{ir}(\fatx_{i,\cdot}).
\label{eq:splittingerrorsecondmoment}
\end{align}
\noindent
While this expression is more complicated than \eqref{eq:splittingerrormean}, the amount of additional work required to compute it is not that great since the complexity of the extra terms are all $\mathcal{O}(N)$ and most of the expensive propensity function evaluations can be  overlapped with the computations involved in \eqref{eq:splittingerrormean}. A more commonly used second order statistic is the variance $V[X] = E[(X-E[X])^2]$. While this quantity is not readily obtained from \eqref{eq:weakerror}, an approximate formula for $V[X]$ can be found from   the observation 
\begin{align}
\nonumber
\Delta V[X] = \Delta E[X^2] - E[X]^2+E_f[X]^2 =\\ 
\nonumber
\Delta E[X^2] + (\Delta E[X])^2 - 2E[X]\Delta E[X] = \\
\Delta E[X^2] - 2x\Delta E[X] + \mathcal{O}(\Delta t^3),
\end{align}
\noindent
where $E_f[X]$ is the expected value under the operator split PDF. Hence   
\begin{align}
\nonumber
\Delta V[X_{is}] =&\quad 0.5\Delta t^2(|D|R(x))_{is}  - \Delta t^2\sum_r n_{irs}\sum_{j \neq i} (d_{ijs}x_{is}\Delta a_{irs}^- - d_{jis}x_{js}\Delta a_{irs}^+)\\
& - 0.5\Delta t^2\sum_{r} n_{irs}^2\sum_{s'=1}^S\sigma_{irs'} - \Delta t^2\sum_{j \neq i} d_{ijs}\sum_r n_{irs}^2a_{ir}(\fatx_{i,\cdot}) + \mathcal{O}(\Delta t^3).
\end{align}
\noindent 
In the asymptotic regime, $\mathcal{O}(\Delta t^3)$ terms are small by definition and $\Delta V[X_{is}]$ is a good approximation to the true local error in variance.
  
We point out that \eqref{eq:eserrpdf}, and hence \eqref{eq:splittingerrormean} and \eqref{eq:splittingerrorsecondmoment}  are only good estimates of the local error if the $\mathcal{O}(\Delta t^3)$ terms are small. It is well known that an error estimate based on \eqref{eq:eserrpdf} will deteriorate in quality for large timesteps if  $\opD$ or $\opM$ are stiff. As the norm of $\opD$ increases rapidly with finer mesh resolution, this will lead to overly conservative estimates for large  $\Delta t>>h^2/\gamma_s$ (in the non-assymptotic regime). 
For this reason, error estimates that perform better in the non-asymptotic regime have been devised in the PDE case \cite{DescombesDLM07}. However, it is not clear how those approaches would apply to the RDME. In the following section we will demonstrate that \eqref{eq:eserrpdf} is simple enough to lead to a strategy that is both computable in practice during the coarse of a single realization of RDME and local in space so that it has potential to be efficiently implemented in parallel. 

\begin{figure}[H]
\includegraphics[width=15cm]{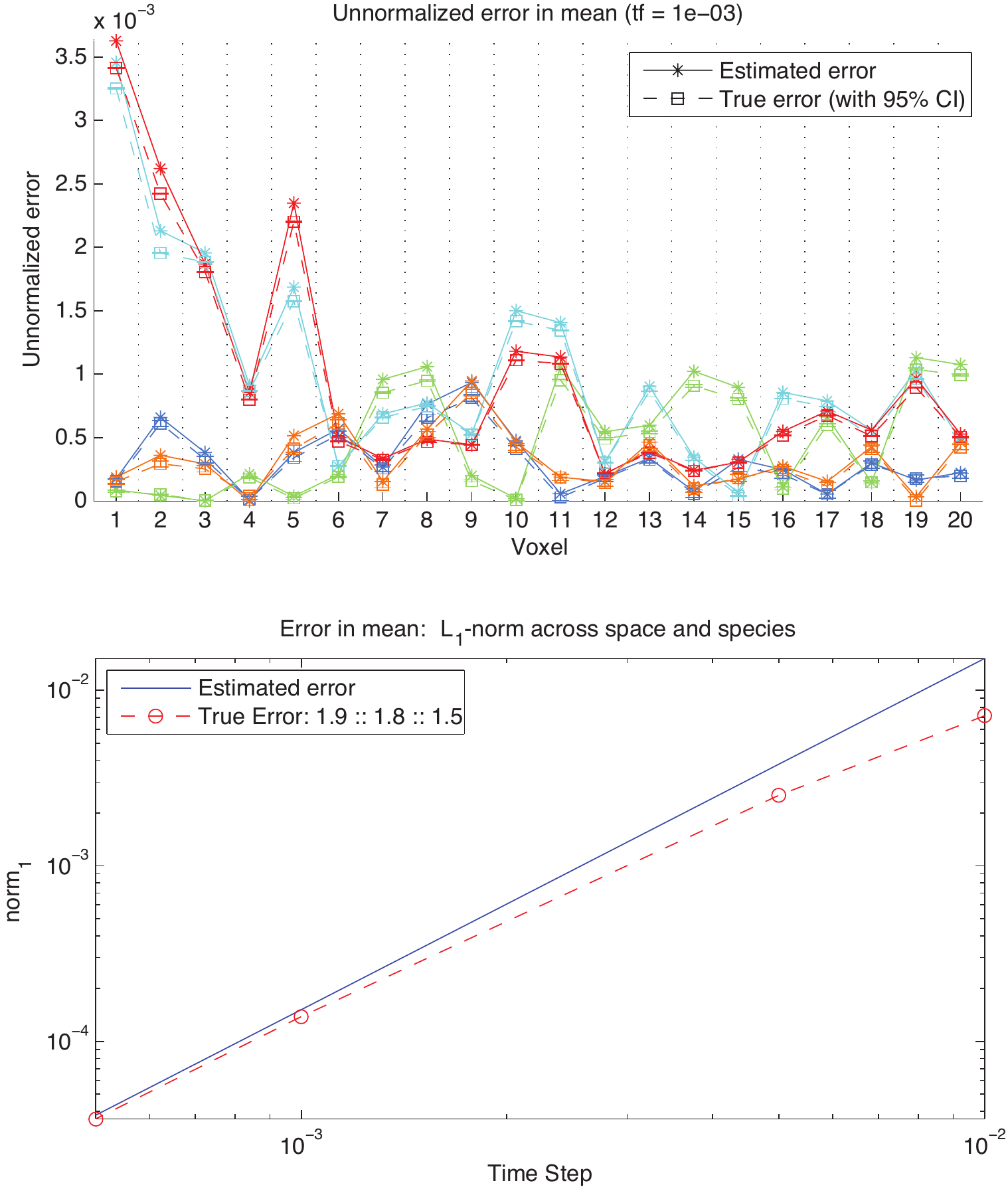} 
\caption{The local error in mean for a single input point x0 vs the estimated error.
{\bf Top:} True (square) and estimated error (*) for each species and voxel for t = 0.001s (voxels are divided by black dotted lines). 
Blue: species 1, cyan: species 2, green: species 3, red: species 4, orange: species 5. 
{\bf Bottom:} The $L_1$-norm across species and voxels of the value in Equation \eqref{eq:localerrormeanl1} for various time-steps. 
The solid blue is the estimated error from Equation \eqref{eq:splittingerrormean}. 
The dashed red line is the true local error. 
The values in the legend are the slope of the curves in the loglog plot. Note that as the time step decreases the estimator converges toward the true error, which in turn approaches the expected $\mathcal{O}(\Delta t^2)$ convergence rate.  
}
\label{fig:fmerrorcomp}
\end{figure}

\begin{figure}[H]
\includegraphics[width=15cm]{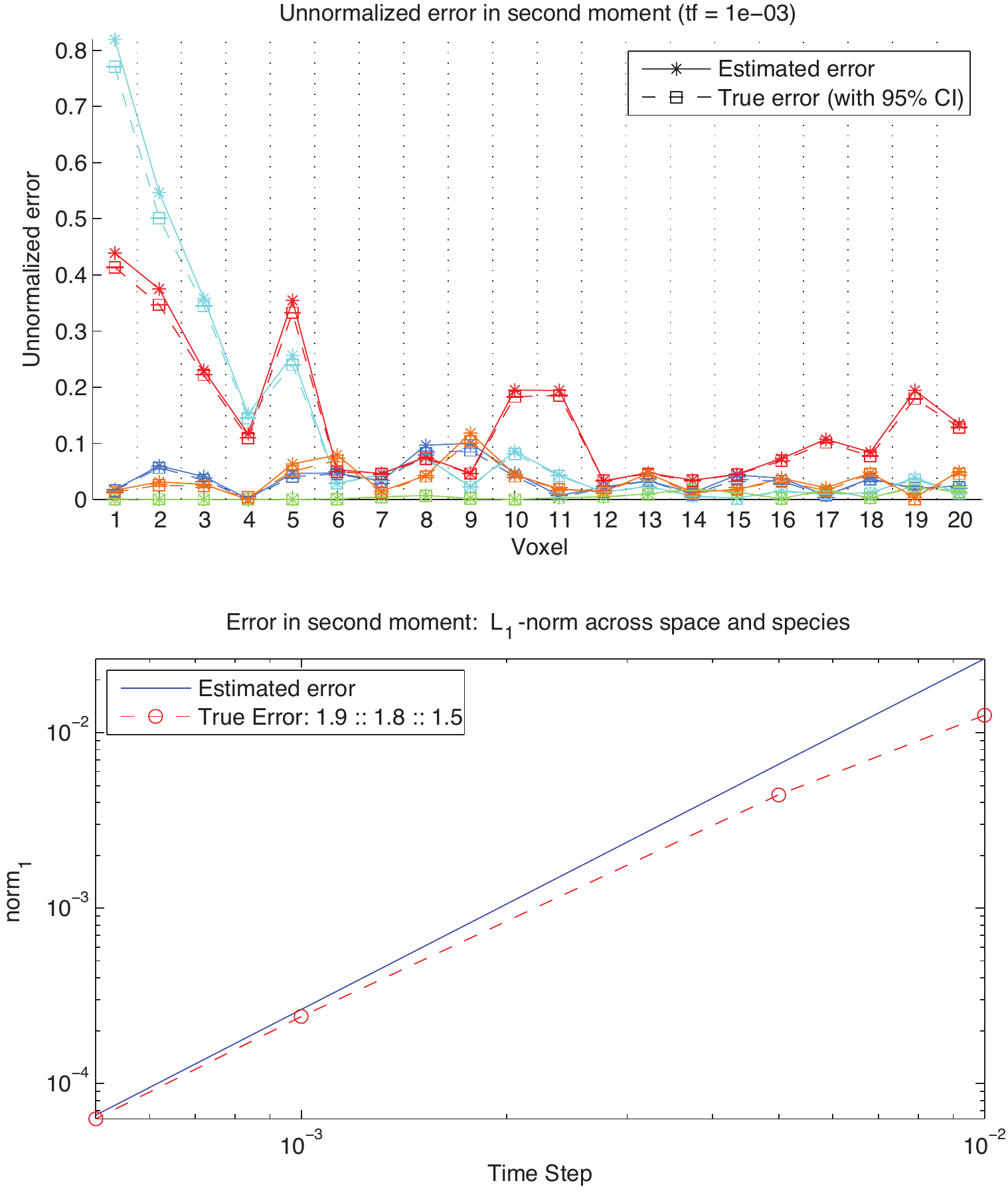} 
\caption{The local error in second moment for a single input point x0 vs the estimated error.
{\bf Top:} True (squares) and estimated error (*) for each species and voxel for t = 0.001s (voxels are divided by black dotted lines). 
Blue: species 1, cyan: species 2, green: species 3, red: species 4, orange: species 5. 
{\bf Bottom:} The $L_1$-norm across species and voxels of the value in Equation \eqref{eq:localerrormeanl1} for various time-steps. 
The solid blue is the estimated error from Equation \eqref{eq:splittingerrorsecondmoment}. 
The dashed red line is the true error. }
\label{fig:smerrorcomp}
\end{figure}

To illustrate the accuracy and correctness of the estimate of the local splitting errors \eqref{eq:splittingerrormean} and \eqref{eq:splittingerrorsecondmoment}, we simulated a model of Min oscillations in \emph{E. Coli}  \cite{BISTAB} in one spatial dimension for a single timestep. To isolate the splitting error, the reaction and diffusion steps were solved exactly by simulating them with SSA and NSM respectively. Fig. \ref{fig:fmerrorcomp} shows the estimated conditional error in mean \eqref{eq:splittingerrormean} and the true conditional error in mean as a function of the spatial coordinate (left) and for different timesteps (right). Fig. \ref{fig:smerrorcomp} shows the corresponding estimate of error in second moment based on \eqref{eq:splittingerrorsecondmoment}. As can be seen, the estimated error accurately captures the true local error, and the quality of the estimate improves for smaller timesteps as expected from the $\mathcal{O}(\Delta t^3)$ error in the estimates. When computing the error in mean, the observed error is a combination of sampling error caused by a finite number of realizations, and the error caused by operator splitting. For small splitting errors, a large number of realizations are necessary to distinguish the sampling error from the splitting error, especially if we want to measure the error in e.g. the $L_1$ norm, in which case the variances in the different voxels add up. 
We note that to achieve tight confidence intervals on the true error for this simple system (see to the top row of Figs. \ref{fig:fmerrorcomp} and \ref{fig:smerrorcomp}) required $10^{11}$ realizations.

\section{Example: Adaptive Diffusive Finite State Projection method}

To illustrate the use of the local error estimate, we extended a previously published approximate method for the RDME, the Diffusive Finite State Projection (DFSP) method \cite{dfsp} with temporal adaptivity. This results in increased robustness of the solver. 

DFSP relies on Lie-Trotter splitting to separate the reaction and diffusion updates over a discrete timestep $\Delta t$. The reactions are then updated in each voxel using Gillespie's direct method \cite{SSA}. Diffusion is updated by sampling from probability density functions that are precomputed by solving the diffusion master equation \eqref{eq:RDMEdf} locally up to the given timestep $\Delta t$. For an unstructured mesh, there is one such PDF for each vertex in the mesh and for each $\Delta t$. The (spatial) locality of the PDFs are enforced by applying an absorbing boundary condition at a certain distance away from the vertices. The timestep is assumed to be small enough that the majority of probability is located close to the vertex where the molecule is located at the beginning of the timestep so that the error due to the truncation of the statespace is small. It is shown in \cite{dfsp} that this strategy can speed up simulations by effectively aggregating the effects of many fast diffusive transfers over the splitting timestep, and in \cite{URDME_BMC} it is discussed under what conditions one can expect DFSP to be more efficient than NSM for the MinCDE model \cite{Huang:2003}. 

DFSP has been implemented previously as an add-on solver in the URDME framework  \cite{URDME_BMC}. URDME is a modular framework that uses unstructured meshes for spatial stochastic simulations. URDME includes interfaces for handling of the geometry, mesh and computation of diffusion jump rates for the unstructured mesh. It has a modular design which facilitates the implementation of new algorithms as add-on solvers. In the current implementation of DFSP in URDME, $\Delta t$ has to be chosen by trial and error by manually picking a timestep that results in satisfactory error in the computed solutions, judged by \emph{a posteriori} checks.  
 
\begin{figure}[htb]
\centering
\includegraphics[trim=3.25cm 8.25cm 3.cm 8.25cm, clip=true,width=.85\linewidth]{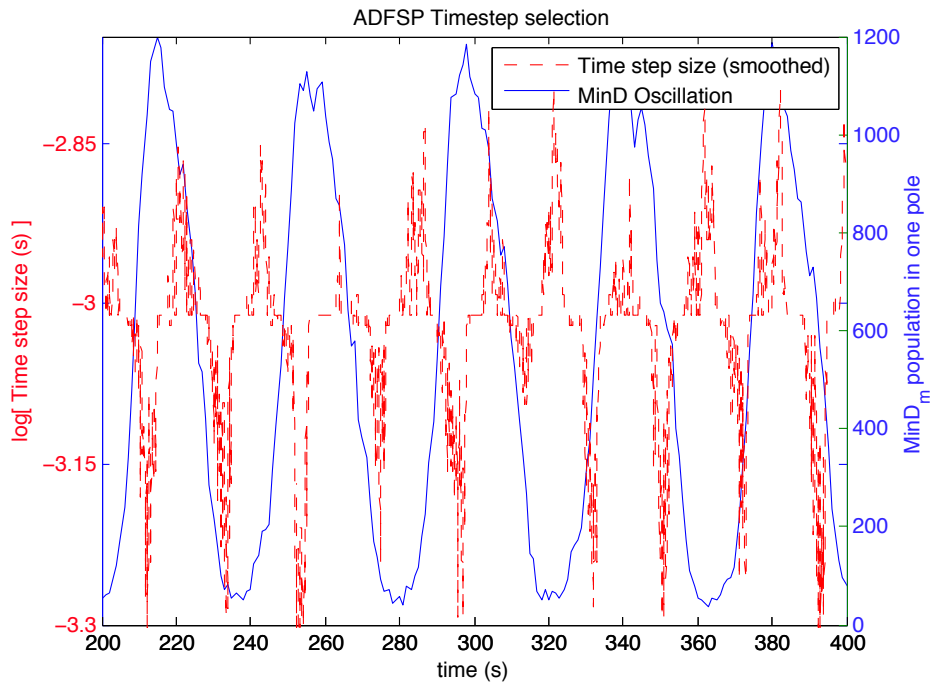}
\caption{The figure shows the time steps selected by the adaptive DFSP method, along with the oscillation pattern of the membrane bound MinD protein for a representative trajectory. Note that the timestep adapts to the dynamics of the MinD oscillation. Due to the oscillatory nature of this model, picking the most conservative timesteps based on the maximal errors would lead the solver to take unnecessarily many steps, resulting in suboptimal performance.  
}
\label{fig:adfsp_timestep_selection}
\end{figure}

We implemented our error estimation strategy in URDME for arbitrary processes (i.e. it uses the same input files as all the core solvers) for structured and unstructured meshes, and supplemented the DFSP solver with an adaptive selection of the timestep. To compute the next proposed timestep, we control the per species error in mean \eqref{eq:splittingerrormean} in the $L_1$ norm such that 
\begin{equation}
\frac{\sum_{i=1}^N |\voxel_i| |\Delta E[X_{is}]|}{\sum_{i=1}^N |\voxel_i|x_{is}} \le \epsilon_s, \quad  s=1,\ldots, S
\label{eq:localerrormeanl1}
\end{equation}
\noindent
with $\epsilon_s$ being a user specified relative error tolerance. The $L_1$ norm is a natural choice in the context of DFSP since the error in the diffusion lookup-tables has a natural bound in the $l_1$ norm \cite{dfsp}. Figure \ref{fig:adfsp_timestep_selection} shows an overlay of the pole-to-pole oscillation pattern of membrane bound MinD in {E. coli} along with time steps selected by the adaptive step size selection algorithm. The timesteps are themselves stochastic variables and fluctuate during the course of the simulation. In the figure they have been smoothed by a windowed average to more clearly visualize how they adapt to the dynamics of the solution. 

We simulated the MinCDE model in 3D in three different ways: using the adaptive timestep control, using a fixed manually determined timestep and using the exact NSM solver, and compared the execution times.  In DFSP, the diffusion step is conducted by computing new positions for each molecule individually by sampling from a precomputed lookup-table. All the molecules' new positions can then be sampled in parallel. During the reaction update, all the voxels can be treated independently and in parallel.  Note that the local error estimators used here have the same characteristics when it comes to the spatial access pattern as  the DFSP diffusion operator. Thus the overall adaptive algorithm should parallelize quite well on a shared memory multicore architecture. To demonstrate this, we implemented a na{\"i}ve parallel version using openMP. Table \ref{tab:dfsp_performance} shows simulation speeds for the adaptive DFSP method, with error estimation conducted every tenth timestep, on a machine with an Intel(R) Core(TM) i7-2600 CPU @ 3.40GHz (Nehalem) processor  and 6GB of RAM using one core and using 4 cores (8 threads). For comparison, we also include simulation times of the (serial) NSM solver in URDME. As can be seen, the DFSP algorithm shows an almost ideal speedup on 4 cores, and for this error tolerance, the adaptive code is roughly two times faster than NSM on a single core.  

\begin{table}[htp]
\centering
\begin{tabular}{c|cc|c}
&\multicolumn{2}{|c|}{Adaptive DFSP ($\epsilon_s = 0.05$)}  & NSM \\ 
\# Voxels & 1 core & 4 cores & (serial)\\ \hline
588 & 124 & 35 &  337 \\
1009  &  228 & 65 & 537\\
2818 & 757 &190 & 1282
\end{tabular}
\caption{Simulation time(s) for the adaptive version of the DFSP solver and for the exact NSM solver for the MinCDE problem with varying mesh resolution. The trajectories were simulated to a final time of $2000$ s, and the state was sampled every second.}
\label{tab:dfsp_performance}
\end{table}

\section{Discussion}

In this work we have derived local error estimators for first order operator splitting in stochastic reaction-diffusion simulations based on the RDME. 
Operator splitting provides a way to decouple reactions and diffusion in spatial stochastic simulations, and the local error estimates enable this to be accomplished with a controlled error. 
The advantage of decoupling the reaction and diffusion operators for simulation are two fold. 
The first is that decoupling the operators allows for approximate methods to be applied to the diffusion operator to reduce the cost of frequent diffusive transfers. The method we considered as an example in this paper falls into this category. Enabling temporal adaptivity for such approximate and hybrid methods has several benefits. Possibly the most important is the robustness it adds to the solver. 
While it is possible to prescribe an error tolerance that is likely to give reasonable result across a range of models, the same is not possible by selecting a fixed timestep \emph{a priori}. 
Adaptivity has long been available in state-of-the art software for e.g. numerical solution of ODEs and PDEs, but was previously lacking for this class of approximate methods for the RDME. The other major benefit is one of efficiency; using a globally conservative fixed timestep will lead to sub-optimal simulation speeds.  We found that the adaptive version of the DFSP method presented here was roughly 1.5 times faster than the corresponding fixed timestep method using a conservative timestep for the oscillatory MinCDE model problem we considered. Even though this experiment is somewhat artificial since the conservative fixed timestep could not be known before running the adaptive code (which again illustrates the utility of adaptivity), it illustrates the efficiency benefits of an adaptive solver.  

The second advantage of operator splitting in the context of the RDME is that splitting the operators converts the largely serial exact simulation method to one with better potential for parallel implementations. 
We have here demonstrated that the local error estimates can be efficiently applied also in a parallel implementation of an approximate method for the RDME. The parallelelization aspects of the current work will be explored in greater detail in a forthcoming publication.  

Finally, operator splitting allows the flexibility of coupling different types of models (possibly at different scales) and solvers. Thus this error estimator opens avenues of research in the development of hybrid algorithms with controlled errors.

\section{Acknowledgment}

Per L{\"o}tsedt provided valuable comments on the manuscript. This work was funded by National Science Foundation (NSF) Award No. DMS-1001012, ICB Award No. W911NF-09-0001 from the U.S. Army Research Office, NIBIB of the NIH under Award No. R01-EB014877-01, and (U.S.) Department of Energy (DOE) Award No. DE-SC0008975. The content of this paper is solely the responsibility of the authors and does not necessarily represent the official views of these agencies.

\bibliographystyle{abbrv}
\bibliography{opsplitting_archive.bib}

\end{document}